\documentclass[fleqn,12pt]{article}
\usepackage[latin1]{inputenc}
\usepackage{graphicx,amsmath,amsfonts,amssymb,amscd,bezier,amstext,makeidx}

\newtheorem{theorem}{Special Theorem}[section]
\newtheorem{defi}[theorem]{Definition}
\newtheorem{theo}[theorem]{Theorem}
\newtheorem{lemm}[theorem]{Lemma}
\newtheorem{corol}[theorem]{Corollary}
\newtheorem{remark}[theorem]{Remark}

\newtheorem{prop}[theorem]{Proposition}

\title{Best simultaneous Diophantine approximations for some pairs of algebraic numbers}
\author{Gustavo A. Pavani }
\date{August 26, 2021}

\begin{document}

\maketitle

\begin{abstract}
\noindent We compute the sequence of best approximations for the vector $(1/\beta,1/\beta^{2})$ where $\beta$ is the real dominant root of the polynomial $P_{a,b}(x)=x^{3}-ax^{2}-bx-1$, where $a$ and $b$ are integers satisfying $-a+1 \leq b \leq -2$. In this case $\beta$ is a Pisot number which does not satisfy the finiteness property. 
\end{abstract}

\section{Introduction}

Our aim is to find the sequence of best approximations for pairs of Pisot numbers which does not satisfy a combinatorial property called Finiteness Property. To achieve this we will use the topological properties of a class of Rauzy fractals related to these type of Pisot numbers.\\

Let $v$ be an element of  $\mathbb{R}^{d}$, $d \geq 2$, $\mathcal{N}$ be a norm in $\mathbb{R}^{d}$ and $q$ be an integer number. Let us define

\centerline{$\mathcal{N}_{0}(qv)=min\{\mathcal{N}(qv-(p_{1},...,p_{d})), (p_{1},...,p_{d})\in\mathbb{Z}^{d}\}$.}

Let $(q_{n})_{n \geq 0}$ be an increasing sequence of integer numbers. We say that $(q_{n})_{n\geq 0}$ is the sequence of best approximations of the vector $v$ for the norm $\mathcal{N}$, if for all $n\in \mathbb{Z}^{+}$ and $0<q<q_{n}$, one have $\mathcal{N}_{0}(q_{n}v)<\mathcal{N}_{0}(qv)$.

This sequence depends on the norm $\mathcal{N}$ and there is no general algorithm that furnishes the sequence of best approximations for all elements of $\mathbb{R}^{d}$. When $d=1$ the problem is already solved: the sequence of best approximations is given by the classical algorithm of one-dimensional continued fractions. But, when $n\geq2$ the problem of finding sequence of best approximations becomes a very difficulty one because there is no algorithm that provides this sequence for all vectors (see \cite{Lagarias}). We can consider this question for classes of vectors, in particular, for vectors which coordinates are cubic Pisot numbers. A Pisot number is an algebraic integer greater than 1 such that its conjugates have modulus less than 1.

The first results is this sense was obtained in \cite{Chekhova}, \cite{Chevallier}, and \cite{HubMess}. These results were obtained by using the topological properties of the so-called Rauzy fractals (figures 1 and 2), which can be defined by means of $\beta-$representations (a sequence of integer numbers) \cite{Rauzy1982}. If a $\beta$-representation ends up with infinitely many zeros, it is said to be finite. We say that a Pisot number $\beta$ has the Finiteness Property (or Property (F) ) if the $\beta-$representation of every nonnegative element of $\mathbb{Z}[\beta]$ is finite. 

The works mentioned before only considered Pisot numbers having the Property (F). In this paper we will prove the following:

\begin{theo} \label{bestseq} Let $\beta$ be a cubic Pisot unit, with complex conjugates satisfying the equation $X^{3}=aX^{2}+bX+1$, $-a+1 \leq b \leq -2$. Let $(T_{n})_{n \geq 0}$ be the recurrent sequence of integer numbers defined by

\begin{center} $T_{0}=1$, $T_{1}=a$, $T_{2}=a^{2}+b, T_{n+3}=aT_{n+2}+bT_{n+1}+T_{n}$, $\forall n \geq 0$.
\end{center}

If $\beta$ does not have the Property (F) then there exists a norm $\mathcal{N}$ in $\mathbb{R}^{2}$ (called Rauzy Norm) and $n_{0} \in \mathbb{Z}^{+}$ such that $(T_{n})_{n \geq n_{0}}$ is the sequence of best approximations for the vector $(1/\beta,1/\beta^{2})$ for the norm $\mathcal{N}$.
\end{theo}

Our theorem is an extension of that one by Messaoudi and Hubert (\cite{HubMess}) for a class of vectors which coordinates are Pisot numbers not having the  Property (F). Hereafter, we will follow the notations of Rauzy, Messaoudi and Hubert.

\section{Numeration systems and Rauzy fractals}

\subsection{ $\beta-$numeration}

Given a real number $\beta>1$, a $\beta$-representation (or $\beta$-expansion) of a number $x \in \mathbb{R}^{+}$ is an infinite sequence $(x_{i})_{i\leq k}$, where $k \in \mathbb{Z}$, $x_{i} \geq 0$ such that $x=\sum_{i=-\infty}^{k}x_{i}\beta^{i}$. The digits $x_{i}$  can be computed using the greedy algorithm as follows (see \cite{Parry, Frougny2000} for details): denote by $\left\lfloor x\right\rfloor$ and $\left\{x\right\}$ the integer and fractional parts of the number $x$. There exists $k \in \mathbb{Z}$ such that $\beta^{k} \leq x < \beta^{k+1}$. Let $x_{k}=\left\lfloor x/\beta^{k}\right\rfloor$ and $q_{k}=\left\{x/\beta^{k}\right\}$. Then, for $i<k$, put $x_{i}=\left\lfloor \beta q_{i+1}\right\rfloor$ and $q_{i}=\left\{\beta q_{i+1}\right\}$. We obtain $x=x_{k}\beta^{k}+x_{k-1}\beta^{k-1}+\cdots$. If $k<0$ ($x<1$) we put $x_{0}=x_{1}=\cdots=x_{k+1}=0$. As we said it before, if a $\beta$-representation ends up with infinitely many zeros, it is said to be finite and the ending zeros can be omitted. Then, the sequence will be denoted by $(x_{i})_{n\leq i \leq k}$ or $x_{k}\cdots x_{n}$. The digits $x_{i}$ belong to the set $B=\{0, \cdots, \beta\}$, if $\beta$ is an integer, or to the set $B=\{0, \cdots, \left\lfloor \beta\right\rfloor\}$, otherwise. 

Cubic Pisot units were classified by Akiyama in \cite{Akiyama2000} as being exactly the set of dominant roots of the polynomial $P_{a,b}(x)=x^3-ax^{2}-bx-1$, satisfying one of the following conditions

\hspace{0.50cm} \textbf{i)} $1 \leq b \leq a$ and $d(1,\beta)=.ab1$;

\hspace{0.50cm} \textbf{ii)} $b=-1$, $a \geq 2$ and $d(1,\beta)=.(a-1)(a-1)01$;

\hspace{0.50cm} \textbf{iii)} $b=a+1$ and $d(1,\beta)=.(a+1)00a1$;

\hspace{0.50cm} \textbf{iv)} $-a+1 \leq b \leq -2$ and $d(1,\beta)=.(a-1)(a+b-1)(a+b)^{\infty}$,

where $(a+b)^{\infty}$ represents the periodic expansion $(a+b)(a+b)(a+b)\ldots,$ and $d(1,\beta)$ is the R\'{e}nyi $\beta$-representation of 1 (see \cite{Renyi} for the definition). \\

Let Fin($\beta$) be the set of nonnegative real numbers that have a finite $\beta$-representation. We say that a Pisot number $\beta$ has the Finiteness Property (or Property (F) ) if $\mathbb{Z}[\beta] \cap [0,+\infty[ \subset$ Fin($\beta$). Therefore, the Pisot numbers in the sets \textbf{i)}, \textbf{ii)} and \textbf{iii)} have the Property (F), while the Pisot numbers in \textbf{iv)} have not.\\

Let us suppose that $\beta$ is a cubic Pisot unit which does not satisfy the Property (F) and let us denote by $\alpha$ and $\lambda$ its Galois conjugates. Let $P_{a,b}(x)=x^{3}-ax^{2}-bx-1$ be the minimal polynomial of $\beta$. We will introduce a generalization of numeration systems induced by the $\beta$-expansions applied on integer numbers.\\

Let $(T_{n})_{n\geq0}$ be the recurrent sequence defined by $T_{0}=1$, $T_{1}=a$, $T_{2}=a^{2}+b$, ${T_{n+3}=aT_{n+2}+bT_{n+1}+T_{n}}$, satisfying the condition $-a+1 \leq b \leq -2$ for all $ n \geq 0$. We have the following lemma:

\begin{lemm} The sequence $(T_{n})_{n \geq 4}$ satisfies

\begin{center} $T_{n}=(a-1)T_{n-1}+ (a+b-1)T_{n-2}+(a+b)T_{n-3}+\cdots +(a+b)T_{1}+(a+b+1)T_{0}$,  \end{center}

for all $n \geq 4$. \end{lemm}

\noindent \textbf{Proof.} The proof is by recurrence on $n$. $\Box$\\

\begin{prop} \label{todoint} Every nonnegative integer $n$ can be uniquely expressed as $n=\sum_{i=0}^{N}d_{i}T_{i}$, where $d_{i} \in \{0, \ldots, a-1\}$ and $d_{j}d_{j-1}\cdots d_{j-k} \leq_{lex}(a-1)(a+b-1)(a+b)\cdots (a+b)$, for all $j \geq k \geq 0$, where ``$\leq_{lex}$'' is the lexicographical order. \end{prop}

\textbf{Proof.} We can use the greedy algorithm to obtain the digits $(d_{j})_{0 \leq j \leq N}$. Notice that  $(T_{n})_{n \geq 0}$ is an increasing sequence of natural integers, since $-a+1 \leq b \leq -2$. Hence, by the definition of the greedy algorithm, we can prove that \\

\centerline{$\displaystyle \sum_{i=0}^{j}d_{i}T_{i}<T_{j+1}$,}

\bigskip

for all $0 \leq j \leq N$  (see \cite{Parry}). Thus, 

\begin{equation} d_{j}d_{j-1} \cdots d_{j-k} <_{lex} (a-1)(a+b-1)(a+b) \cdots (a+b)(a+b+1), \forall j \geq k \geq 0. 
\end{equation}

Therefore we obtain that $d_{j} \cdots d_{j-k}\leq_{lex}(a-1)(a+b-1)(a+b) \cdots (a+b)$. $\Box$\\

Let $\frak{D}=\{(d_{i})_{i\geq k}, k\in \mathbb{Z}, \forall n\geq k, d_{n}\cdots d_{n-k}\leq_{lex}(a-1)(a+b-1)(a+b)\cdots (a+b)\}$. By definition, the Rauzy fractal is the set \\

\centerline{$\displaystyle \mathcal{R}:=\mathcal{R}_{a,b}=\left\{\sum_{i=2}^{+\infty}d_{i}\theta_{i}, \, (d_{n})_{n\in \mathbb{Z}} \in \frak{D}\right\}$}

\bigskip

where $\theta_{i}=\alpha^{i}$, if $\alpha \in \mathbb{C}\setminus \mathbb{R}$ or $\theta_{i}=(\alpha^{i},\lambda^{i})$, if $\alpha \in \mathbb{R}$. Notice that  $\mathcal{R} \subset \mathbb{C}$ or $\mathcal{R} \subset \mathbb{R}^{2}$.\\

\textbf{Alternative definition of the Rauzy fractal.} 

Let $N \in \mathbb{Z}^{+}$ and $(d_{j})_{k(N) \geq j \geq 0}$ be a $T$-representation of $N$, i.e., $N=\sum_{i=0}^{k(N)}d_{i}T_{i}$, where $(d_{i})_{i \geq 0} \in \frak{D}$. Let

\centerline{$\delta(N)=N(1/\beta,1/\beta^{2})-(P_{N},Q_{N})$,}

where

\centerline{$P_{N}=\displaystyle \sum_{j=1}^{k(N)}d_{j}T_{j-1}, \,\, Q_{N}=\displaystyle \sum_{j=2}^{k(N)}d_{j}T_{j-2}$.}

Let us consider the matrix $B$ defined as\\

\centerline{$B=\left(
                 \begin{array}{cc}
                   -b/\beta & -1/\beta \\
                   1-b/\beta^{2} & -1/\beta^{2} \\
                 \end{array}
               \right)$
.}

We have the following property.\\

\begin{lemm}
For all $n \geq 2$,\\

\centerline{$B\left(
                 \begin{array}{c}
                   T_{n}/\beta-T_{n-1} \\
                   T_{n}/\beta^{2}-T_{n-2} \\
                 \end{array}
               \right)=\left(
                         \begin{array}{c}
                           T_{n+1}/\beta-T_{n} \\
                           T_{n+1}/\beta^{2}-T_{n-1} \\
                         \end{array}
                       \right)$
}

\noindent where $(T_{n})_{n\geq0}$ is the sequence defined previously.
\end{lemm}

\textbf{Proof.} The proof is not difficult and it is made by induction. $\Box$

\begin{corol}
If $N=\sum_{j=0}^{k(N)}d_{j}T_{j}$ then $\delta(N)=\sum_{j=0}^{k(N)}d_{j}B^{j}\delta(1)$, where $\delta(1)=(1/\beta,1/\beta^{2})$.
\end{corol}

Thus, the Rauzy fractal is the set\\

\centerline{$\mathcal{E}$ $=\overline{\{\delta(N);\, N \in \mathbb{Z}^{+}\}} \subset \mathbb{R}^{2}$.}

\bigskip
This set was introduced by G. Rauzy in 1982 (\cite{Rauzy1982}). Since then, this set and its generalizations have been extensively studied due to its strong connections with many fields of mathematics such as Dynamical Systems and Number Theory.

As mentioned before, our results highly depends on the topological properties of the class of Rauzy fractals associated to the Pisot numbers not having the Property (F). The topological and arithmetical properties of this class of Rauzy fractals were studied in details in \cite{PaperGustavo} (see also \cite{TeseGustavo}). In particular, it was proved that $\mathcal{R}_{a,b}$ has at least $6+2(K-1)$ neighbors, where $\displaystyle K=\left[\frac{a-1}{a+b+1}\right]$. Let us remind that an element $u$ of a lattice $\Lambda$ is a neighbor of $\mathcal{R}$ if, and only if, $\mathcal{R} \cap \mathcal{R}+u \neq \emptyset$. We will assume that we are in the case of 6 neighbors. \\

\begin{figure}[h!]
\centering

\hspace{2.0cm}\includegraphics[scale=0.15]{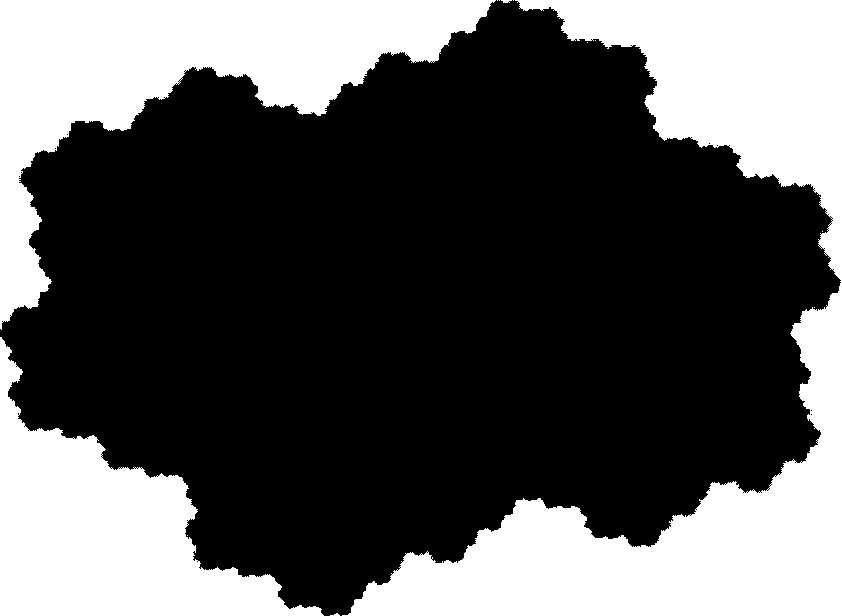} 
\hspace{3.0cm}
 \includegraphics[scale=.3]{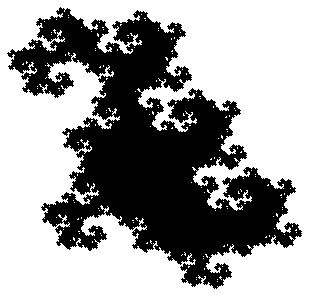} 
\caption{The set $\mathcal{R}_{1,1}$ (Classic Rauzy Fractal). \hspace{1.0cm} Figure 2: The set $\mathcal{R}_{3,-2}$.}
\end{figure}

\subsection{Rauzy Norm}

Let  $M=\left(
          \begin{array}{cc}
            \lambda+b/\beta & 1/\beta \\
            -\alpha-b/\beta & -1/\beta \\
          \end{array}
        \right)$. One can check that the matrix $B$ defined before is similar to the matrix $\left(
                                   \begin{array}{cc}
                                     \alpha & 0 \\
                                     0 & \lambda \\
                                   \end{array}
                                 \right)$ and it satisfies
																
\centerline{$MB=\left(
                  \begin{array}{cc}
                    \alpha & 0 \\
                    0 & \lambda \\
                  \end{array}
                \right)M.$
}

We will consider the case when $\beta$ has complex conjugates. In this case, the Rauzy Norm  $\mathcal{N}$ is defined by:

\centerline{$\mathcal{N}(x)=|(\alpha+b/\beta)x_{1}+x_{2}/\beta|, \forall x=(x_{1},x_{2})\in \mathbb{R}^{2}$.}

\begin{remark}
\noindent Notice that $\mathcal{N}(x)=|\pi_{1}Mx|$, where $\pi_{1}(x,y)=x$, for all $(x,y)\in \mathbb{R}^{2}$, and it is not difficult to verify that $\mathcal{N}(Bx)=|\alpha|\mathcal{N}(x)$, for all $x \in \mathbb{R}^{2}$.
\end{remark}

\begin{corol}
Let $(d_{n})_{n \geq 0}$ be a $T-$representation. Then 

\centerline{$\mathcal{N}(\sum_{n=0}^{+\infty}d_{n}B^{n}\delta(1))=\mathcal{N}(\delta(1)|\sum_{n=0}^{+\infty}d_{n}\alpha^{n}|$.}

In particular, 

\centerline{$\mathcal{N}(\delta(T_{n}))=\mathcal{N}(\delta(1))|\alpha^{n}|$, for all $n \geq 0$.}

\end{corol}


\begin{remark} All the omitted proofs can be found in \cite{TeseGustavo}.
\end{remark}

\section{Sequence of best approximations}
Let us prove the Theorem \ref{bestseq}. To prove it, it will be necessary several auxiliary results. Let us begin with the theorem:

\begin{theo} \label{teo10}
There exists a real number $c >0$ such that for all $q \in \mathbb{Z}^{+}$ and for all $g \in \mathbb{Z}^{2}$, if $\mathcal{N}(q(1/\beta,1/\beta^{2})-g)<c$ implies that $q(1/\beta,1/\beta^{2})-g=\delta(q)$ or $q(1/\beta,1/\beta^{2})-g=\delta(q)-(1,1)$.
\end{theo}

As a corollary we obtain,

\begin{corol}
There exists a real number $c>0$ such that for all  $q \in \mathbb{Z}^{+}$, if $\mathcal{N}_{0}(q(1/\beta,1/\beta^{2}))<c$ then $\mathcal{N}_{0}(q(1/\beta,1/\beta^{2}))=\mathcal{N}(\delta(q))$.
\end{corol}

The following proposition is an adaptation from that one found in \cite{HubMess}, but we will put it here for a better comprehension.\\

\begin{prop}
There exists a linear and bijective application from $\mathbb{R}^{2}$ to $\mathbb{C}$ such that  $f(\cal{E})=\mathcal{R}$ and $f(\mathbb{Z}^{2})=\mathbb{Z}+\mathbb{Z}\alpha$.
\end{prop}
\noindent \textbf{Proof}. Since $(T_{n})_{n \geq 0}$ is recurrent, there exists $h, \ell \in \mathbb{C}$ such that $T_{n}=h\beta^{n}+\ell\alpha^{n}+\overline{\ell}\overline{\alpha}^{n}$ for all $n \in \mathbb{Z}^{+}$. Therefore, for all $n \in \mathbb{N}$,\\

\centerline{$\delta(T_{n})=\left(\begin{array}{c}
                   T_{n}/\beta-T_{n-1} \\
                   T_{n}/\beta^{2}-T_{n-2} \\
                 \end{array}
               \right)=\left(
                         \begin{array}{c}
                           c\alpha^{n+2}+\overline{c}\overline{\alpha}^{n+2} \\
                           d\alpha^{n+2}+\overline{d}\overline{\alpha}^{n+2} \\
                         \end{array}
                       \right),$}
											
\vspace{0.5cm}
											
where $ c=\frac{\ell}{\alpha^{2}}(1/\beta-1/\alpha)$, and $d=\frac{\ell}{\alpha^{2}}(1/\beta^{2}-1/\alpha^{2})$. Solving the system defined by the initial values of $T_{0}$, $T_{1}$ and $T_{2}$ we obtain that

\centerline{$\displaystyle \ell=\frac{\alpha^{2}}{(\alpha-\overline{\alpha})(\alpha-\beta)}$.}

Thus,

\centerline{$c\overline{d}-d\overline{c}=\displaystyle \frac{\alpha^{2}\overline{\alpha}^{2}}{\overline{\alpha}-\alpha}\neq 0.$}

Let us set $g(z)=\left(      \begin{array}{c}
                           cz+\overline{cz} \\
                           dz+\overline{dz} \\
                         \end{array}
                       \right)$ for all $z \in \mathbb{C}$ and $f=g^{-1}$. We obtain that $f(\mathcal{E})=\mathcal{R}$. On the other hand, we have

\centerline{$g(1)=\left( \begin{array}{c}
                           c+\overline{c} \\
                           d+\overline{d} \\
                         \end{array}
                       \right)=\delta(T_{-2})= \left( \begin{array}{c}
                           T_{-2}/\beta-T_{-3} \\
                           T_{-2}/\beta^{2}-T_{-4} \\
                         \end{array}
                       \right)$.}

Since $T_{-1}=T_{-2}=0$, $T_{-3}=1$ and $T_{-4}=-b$, we obtain that $g(1)= \left( \begin{array}{c}
                           -1 \\
                           b \\
                         \end{array}
                       \right)$, and  

\centerline{$g(\alpha)=\left( \begin{array}{c}
                           c\alpha+\overline{c\alpha} \\
                           d\alpha+\overline{d\alpha} \\
                         \end{array}
                       \right)=\delta(T_{-1})= \left( \begin{array}{c}
                           T_{-1}/\beta-T_{-2} \\
                           T_{-1}/\beta^{2}-T_{-3} \\
                         \end{array}
                       \right)= \left( \begin{array}{c}
                           0 \\
                           -1 \\
                         \end{array}
                       \right)$. }
											
Therefore,  $\mathbb{Z}+\mathbb{Z}\alpha=f(\mathbb{Z}^{2})$. $\Box$\\

The next Lemma is essential for the statements to come. Its proof is in the Annex.

\begin{lemm} \label{rzzalpha} $\mathcal{R} \cap (\mathbb{Z}+\mathbb{Z}\alpha)=\{0,-1-(b+1)\alpha\}$.\end{lemm}


For the next proposition, the interior of a set $X$ will be denoted by $int(X)$ and the distance between two sets $X$ and $Y$ will be denoted by $dist(X,Y).$

\begin{prop} \label{zerointerior} 0 $\in$ int$(\mathcal{R} \cup (\mathcal{R}+1+(b+1)\alpha))$.
\end{prop}

\noindent \textbf{Proof.} Let $c=dist(\mathcal{R},\mathbb{Z}+\mathbb{Z}\alpha \setminus \{0,-1-(b+1)\alpha\})$. Let $\mathcal{B}(0,c)$ be the ball of radius $c$ and center in 0. There exists $ n\in \mathbb{N}, p,q \in \mathbb{Z}$ such that $n\alpha^{2}+p\alpha+q \in \mathcal{B}(0,c)$. Then, $|n\alpha^{2}+p\alpha+q|=|n\alpha^{2}+p_{n}\alpha+q_{n}-((p_{n}-p)\alpha+(q_{n}-q))|<c$. Since $n\alpha^{2}+p_{n}\alpha+q_{n} \in \mathcal{R}$, then $(p_{n}-p)\alpha+(q_{n}-q)\in\{0,1+(b+1)\alpha\}$. We have the following cases:

\textbf{Case 1.} If $(p_{n}-p)\alpha+q_{n}-q=0$ then $p_{n}=p$ and $q_{n}=q$. Thus, $n\alpha^{2}+p\alpha+q \in \mathcal{R}$.

\textbf{Case 2.} If $(p_{n}-p)\alpha+q_{n}-q=1+(b+1)\alpha$ then $p=p_{n}-1$ and $q=q_{n}-(b+1)$. Thus, $n\alpha^{2}+p\alpha+q=n\alpha^{2}+p_{n}\alpha+q_{n}-1-(b+1)\alpha \in \mathcal{R}+(1+(b+1)\alpha)$. $\Box$ \\

\noindent \textbf{Proof of Theorem \ref{teo10}.} We have that $\mathcal{R} \cap \mathbb{Z}+\mathbb{Z}\alpha=\{0,-1-(b+1)\alpha)\}$. Hence, $\mathcal{E} \cap \mathbb{Z}^{2}=(f^{-1}(0),f^{-1}(-1-(b+1)\alpha))=(0,(1,1))$. Then there exists a real number $c>0$ such that for all $g \in \mathbb{Z}^{2}$, inf$_{x \in \mathcal{E}}\mathcal{N}(g-x)<c$ implies that  $g=0$ or $g= (1,1)$. Let us suppose that $\mathcal{N}(q(1/\beta,1\beta^{2})-g)<c$. Since $\delta(q)-q(1/\beta,1/\beta^{2}) \in \mathbb{Z}^{2}$ and $\delta(q) \in \mathcal{E}$, then $\mathcal{N}(\delta(q)-(\delta(q)-q(1/\beta,1/\beta^{2})-g)<c$. Therefore $\delta(q)-q(1/\beta,1/\beta^{2})+g=0$ or $\delta(q)-q(1/\beta,1/\beta^{2})+g=(1,0)$, that is, $\delta(q)=q(1/\beta,1/\beta^{2})-g$ or $q(1/\beta,1/\beta^{2})-g=\delta(q)-(1,1)$. $\Box$\\

Now we are in order to prove the main theorem. \\

\noindent \textbf{Proof of Theorem \ref{bestseq}.} Let us prove that there exists $n_{0} \in \mathbb{N}$ such that $(T_{n})_{n \geq n_{0}}$ is the sequence of best approximation  for the vector $(1/\beta,1/\beta^{2})$ for the norm $\mathcal{N}$.\\

Let $c$ be a real number as defined in Theorem \ref{teo10}. Let $n_{0} \in \mathbb{N}$ such that $\mathcal{N}(\delta(T_{n_{0}}))=\kappa|\alpha^{n_{0}}|<c$, where $\kappa=\mathcal{N}(\delta(1))$. Take $n \geq n_{0}$ and $0 <q <T_{n}$. Let us prove that $\mathcal{N}_{0}(T_{n}(1/\beta,1/\beta^{2}))<\mathcal{N}_{0}(q(1/\beta,1/\beta^{2}))$. Suppose that 

\centerline{$\mathcal{N}_{0}(q(1/\beta,1/\beta^{2})) \leq \mathcal{N}_{0}(T_{n}(1/\beta,1/\beta^{2})) \leq \mathcal{N}(T_{n}(1/\beta,1/\beta^{2}))=\kappa|\alpha^{n}| \leq \kappa|\alpha^{n_{0}}|<c$.}

 Hence, $\mathcal{N}_{0}(q(1/\beta,1/\beta^{2}))<c$. Now, suppose that $\mathcal{N}_{0}(q(1/\beta,1/\beta^{2}))=\mathcal{N}(q(1/\beta,1/\beta^{2})-g)$. Thus, by Theorem \ref{teo10} we obtain that $q(1/\beta,1/\beta^{2})-g=\delta(q)$ or $q(1/\beta,1/\beta^{2})-g=\delta(q)-(1,1)$. Let us treat these cases separately.

\noindent \textbf{Case 1.} $q(1/\beta,1/\beta^{2})-g=\delta(q)$. In this case, 

\begin{equation} \label{eq1} \hspace{3cm}\mathcal{N}_{0}(q(1/\beta,1/\beta^{2}))=\mathcal{N}(\delta(q)) <\mathcal{N}(\delta(T_{n})). \end{equation}

Since $q<T_{n}$, then $q=\sum_{i=0}^{n-1}d_{i}T_{i}$. Thus, 

\begin{equation} \label{eq2} \hspace{4cm}\displaystyle \mathcal{N}(\delta(q))=\kappa\left|\sum_{i=0}^{n-1}d_{i}\alpha^{i}\right|. \end{equation}

From (\ref{eq1}) and (\ref{eq2}) we obtain that $|\sum_{i=0}^{n-1}d_{i}\alpha^{i}|<|\alpha^{n}|$.

On the other hand, $\sum_{i=0}^{n-1}d_{i}\beta^{i} \in \mathbb{Q}$, otherwise we would have $\sum_{i=0}^{n-1}d_{i}\beta^{i}=r/s$, where $r,s \in \mathbb{Z}$. Hence, $N\beta^{2}+p_{N}\beta+q_{N}=r/s$, that is, $sN\beta^{2}+p_{N}s\beta+q_{N}s=r$. Thus, $sN=p_{N}s=0$ and $q_{N}s=r$ and then $s=0$. Absurd. Therefore, $\lambda = \sum_{i=0}^{n-1}d_{i}\beta^{i}$ is an algebraic number of degree 3, since $\mathbb{Q}(\lambda) \subset \mathbb{Q}(\beta)$. Thus $\sum_{i=0}^{n-1}d_{i}\alpha^{i}$, $\sum_{i=0}^{n-1}d_{i}\overline{\alpha}^{i}$ and $\sum_{i=0}^{n-1}d_{i}\beta^{i}$ are Galois conjugates, and hence $|\sum_{i=0}^{n-1}d_{i}\alpha^{i}\sum_{i=0}^{n-1}d_{i}\overline{\alpha}^{i}\sum_{i=0}^{n-1}d_{i}\beta^{i}| \in \mathbb{Z}.$ Therefore, 

\centerline{$\displaystyle \left|\sum_{i=0}^{n-1}d_{i}\alpha^{i}\right|^{2} \geq \frac{1}{\sum_{i=0}^{n-1}d_{i}\beta^{i}} \geq \frac{1}{\beta^{n}}=|\alpha|^{2n}$,}

 which implies that 

\centerline{$\displaystyle \left|\sum_{i=0}^{n-1}d_{i}\alpha^{i}\right| \geq |\alpha|^{n}$.}

Absurd.\\

\noindent \textbf{Case 2.} $q(1/\beta,1/\beta^{2})-g=\delta(q)-(1,1)$. In this case, $\mathcal{N}_{0}(q(1/\beta,1/\beta^{2})-g)=\mathcal{N}(q(1/\beta,1/\beta^{2})-g)=\mathcal{N}(q(1/\beta,1/\beta^{2})-(1,1))$. We have that $\mathcal{N}(q(1/\beta,1/\beta^{2})-g)=\kappa|\sum_{k=2}^{n-1}d_{k}\alpha^{k}|$. Let us suppose that $\mathcal{N}(q(1/\beta,1/\beta^{2})-(1,1))=\kappa|\sum_{i=0}^{n-1}d_{i}\alpha^{i}+1+(b+1)\alpha|$.

Claim: 

\centerline{$\mathcal{N}(q(1/\beta,1/\beta^{2})-(1,1)) \geq \mathcal{N}(\delta(T_{n}))=|\alpha^{n}|$.}

Ideed,

\centerline{$\displaystyle \left|1+(b+1)\alpha+\sum_{i=2}^{n-1}d_{i}\alpha^{i}\right|^{2}\cdot\left|1+(b+1)\beta+\sum_{i=2}^{n-1}d_{i}\beta^{i}\right|\geq 1$.}

Hence,

 \centerline{$\displaystyle \left|1+(b+1)\alpha+\sum_{i=2}^{n-1} d_{i}\alpha^{i}\right|\geq \frac{1}{\left|1+(b+1)\beta+\sum_{i=2}^{n-1} d_{i}\beta^{i}\right|}$.}

\noindent On the other hand, take $n\geq n_{0}$ large enough such that $1+(b+1)\beta+\sum_{i=2}^{n-1}d_{i}\beta^{i}>0$. Hence, $1+(b+1)\beta+\sum_{i=2}^{n-1}d_{i}\beta^{i} \leq \sum_{i=0}^{n-1}d_{i}\beta^{i} \leq \beta^{n}$, because $b+1 \leq -1$. Thus, \\

\centerline{$\displaystyle \frac{1}{1+(b+1)\beta+\sum_{i=2}^{n-1} d_{i}\beta^{i}} \geq \frac{1}{\beta^{n}} =|\alpha|^{2n}$. }

Hence, $1+(b+1)\beta+\sum_{i=2}^{n-1}d_{i}\beta^{i} \geq |\alpha|^{n}$. Absurd.

Therefore, $\mathcal{N}_{0}(T_{n}(1/\beta,1/\beta^{2}))<\mathcal{N}_{0}(q(1/\beta,1/\beta^{2}))$ for all the cases. $\Box$\\

\section{Annex}

\subsection{Neighborhood automaton}

This section is dedicated to prove \\

\textbf{Lemma 3.4 } $\mathcal{R} \cap (\mathbb{Z}+\mathbb{Z}\alpha)=\{0,-1-(b+1)\alpha\}$.

To prove this, we need to construct a finite automaton, which will denote by $\mathcal{H}$, that recognizes points with two $\alpha-$representations, that is, a point that can be written in two different ways. These points belong to the boundary of $\mathcal{R}_{a,b}$. Let us begin with a definition.

\begin{defi} A finite automaton is a triple $(S,A,C)$, where $A$ is the alphabet of the automaton, $S$ is the set of the states and $C$ is a subset of $S \times A \times S$. We say that a sequence $(a_{n})_{n \in \mathbb{N}}$ is recognizable by the automaton $(S,A,C)$ if there exists a sequence $(s_{n}) \in A^{\mathbb{N}}$ such that $(s_{i-1},a_{i},s_{i}) \in C$, for all $i \in \mathbb{N}$\end{defi}

\begin{remark} We can add to the automaton a set $I$, with the initial states and a set $F$ with the ending states. In this work we just need the set $I$. 
\end{remark}


Before constructing the automaton we need the following proposition (see \cite{PaperGustavo}):

\begin{prop} \label{xyteorema} Let $x=\sum_{i=l}^{\infty}a_{i}\alpha^{i}$ and $y=\sum_{i=l}^{\infty}b_{i}\alpha^{i}$, where  $l \in \mathbb{Z}$ and $(a_{i})_{i\geq l}$, $(b_{i})_{i\geq l}$ belong to $\mathcal{L}$. Then $x=y$ if, and only if, the set  $\{x(k)-y(k), k \geq l\}$ is finite, where $x(k)=\alpha^{-k+2}\displaystyle \sum_{i=l}^{k}a_{i}\alpha^{i}$ and $y(k)=\alpha^{-k+2} \sum_{i=l}^{k}b_{i}\alpha^{i}$, $\forall k \geq l$.
\end{prop}

As a consequence of this proposition, we have the following result.

\begin{theo} Let $(a_{i})_{i \geq l}$ and $(b_{i})_{i \geq l}$ two distinct elements of $\mathcal{L}$, then $\sum_{i=l}^{\infty}a_{i}\alpha^{i}=\sum_{i=l}^{\infty}b_{i}\alpha^{i}$ if and only if the sequence $((a_{i},b_{i}))_{i \geq l}$ is recognizable by the automaton $\mathcal{H}$.
\end{theo}

Now, let us construct the automaton $\mathcal{H}$ such that  $n,p \in \mathbb{Z}$ and $(\varepsilon_{i})_{i \geq 2}, (\varepsilon^{'}_{i})_{i \geq 2} \in \mathcal{N}$, $n+p\alpha+\sum_{i=2}^{+\infty}\varepsilon_{i}\alpha^{i}=\sum_{i=2}^{+\infty}\varepsilon^{'}_{i}\alpha^{i}$ if, and only if, $(n,0)(p,0)(\varepsilon_{2},\varepsilon^{'}_{2})\cdots$ is an infinite path in the automaton $\mathcal{H}$.

To make it simple, the idea presented in \cite{TeseGustavo} is: let $p$ and $q$ be two states. The set of edges is the set of $(p,(c,d),q) \in S\times\{0,1,...,a-1\}^{2}\times S$ satisfying $q=\frac{p}{\alpha}+(c-d)\alpha^{2}$. The set of initial states is $\{0\}$. We suppose that we are in the case where $\mathcal{R}$ has 6 neihbors. In this case, it was shown also in \cite{TeseGustavo}, Proposition 4.0.4, that the set of the states is $S=\{0, \pm \alpha^{2}, \pm (\alpha+b\alpha^{2}), \pm (\alpha +(b+1)\alpha^{2}),\pm (1+b\alpha+(a-1)\alpha^{2}), \pm (1+(b+1)\alpha+(a+b)\alpha^{2}), \pm (1+(b+1)\alpha+(a+b+1)\alpha^{2})\}$.

 Let $x=\sum_{i=l}^{+\infty}a_{i}\alpha^{i}$ and $y=\sum_{i=l}^{+\infty}b_{i}\alpha^{i}$, where $a=(a_{i})_{i \geq l}$ and $b=(b_{i})_{i \geq l}$ belong to $\cal{D}$. Suppose that $x=y$ and for all $k \geq l$ we set $S_{k}=S_{k}(a,b)=x(k)-y(k)$. We have,

 \begin{equation} \label{Akmais1} \hspace{4.0cm }\displaystyle S_{k+1}=\frac{S_{k}}{\alpha}+(a_{k+1}-b_{k+1})\alpha^{2}. \end{equation}

Let $t$ be the smallest integer such that $a_{t} \neq b_{t}$. Hence $S_{i}(a,b)=0$ for all $i \in \{l,...,t-1\}$. Suppose that $(a_{t},b_{t})=(1,0)$. Then, $S_{t}=\alpha^{2}$. From (\ref{Akmais1}) we deduce that

 $S_{t+1}=\alpha+(a_{t+1}-b_{t+1})\alpha^{2}$ $=\left\{
                                              \begin{array}{ll}
                                                \alpha+b\alpha^{2}, & \hbox{if } (a_{t+1},b_{t+1})=(b,0) \\
                                                \alpha+(b+1)\alpha^{2}, & \hbox{if} (a_{t+1},b_{t+1})=(b+1,0)
                                              \end{array}
                                            \right.$\\
																						
Then, $(\alpha^{2},(\varepsilon+b,\varepsilon),\alpha+b\alpha^{2})$ is an edge that connects the state $\alpha^{2}$ to state $\alpha+b\alpha^{2}$, and $(\alpha^{2},(\varepsilon+b+1,\varepsilon),\alpha+(b+1)\alpha^{2})$ is and edge that connects the state $\alpha^{2}$ to state $\alpha+(b+1)\alpha^{2}$. Continuing with this process, we obtain an infinite path $(S_{i}, (a_{i},b_{i}),S_{i+1})_{i \geq l}$ beginning in the initial state of the automaton. This path will be denoted by $(a_{i},b_{i})_{i \geq l}$. As the set of states $S$ is finite, we obtain a finite automaton (figure below).

\begin{center}
    \includegraphics[width=17cm,angle=0]{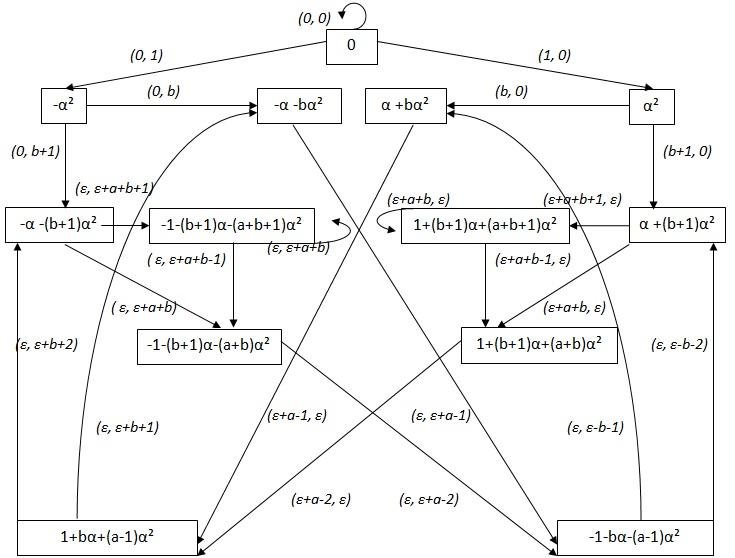}
  \centerline{Figure 3: Neighborhood automaton $\mathcal{H}$.}\vspace{0.5cm}\\

	\end{center}

\textbf{Proof of Lemma \ref{rzzalpha}:} Let us notice that $0$ and $1+b\alpha$ belongs to $\mathcal{R}$. Suppose that there exists $n,p \in \mathbb{Z}$ such that $n+p\alpha=\sum_{i=2}^{\infty}\varepsilon_{i}\alpha^{i}$, where $(\varepsilon_{i})_{i\geq2} \in \mathcal{N}$. Suppose that $n>0$. Then,

\centerline{$(n,0)(p,0)(0,\varepsilon_{2})(0,\varepsilon_{2})\cdots$}

is a path in the automaton $\mathcal{H}$. Absurd, because this sequence is not recognizable by $\mathcal{H}$. \rightline{$\Box$}

\noindent \textbf{Aknowledgments.} I would like to thank Ali Messaoudi for all the help during the writing of this manuscript.

\bibliographystyle{amsplain}
\bibliography{BibSeq}

\small
\textsc{Gustavo Antonio Pavani}

\textsc{Department of Mathematics - State University of Mato Grosso do Sul - UEMS }

\textsc{Address: Rua Walter Hubacher, 138, 79750-000, Nova Andradina, MS, Brazil}

\textsc{E-mail address:} \texttt{gustavo.pavani@uems.br}

\end{document}